\documentclass[11.05pt,a4paper]{amsart}
\usepackage{amssymb,amsmath}
\usepackage{latexsym}
\usepackage{amsthm,amsfonts,amssymb,mathrsfs}
\usepackage{rotating}
\usepackage[leqno]{amsmath}
\usepackage{xspace}
\usepackage[all]{xy}
\usepackage{longtable}

\headsep=1cm \numberwithin{equation}{section}
\newtheorem{theorem}{Theorem}[section]
\newtheorem{lemma}[theorem]{Lemma}

\newtheorem{corollary}[theorem]{Corollary}

\theoremstyle{definition}
\newtheorem{definition}[theorem]{Definition}
\theoremstyle{remark}
\newtheorem{remark}[theorem]{Remark}

\newtheorem{question}[theorem]{Question}

\newtheorem{acknowledgement}{Acknowledgement}

\newcommand{\Ass}{\operatorname{Ass}}
\newcommand{\im}{\operatorname{im}}
\newcommand{\grade}{\operatorname{grade}}
\newcommand{\Min}{\operatorname{Min}}

\newcommand{\Spec}{\operatorname{Spec}}

\newcommand{\cd}{\operatorname{cd}}

\newcommand{\id}{\operatorname{id}}

\newcommand{\Gid}{\operatorname{Gid}}
\newcommand{\pd}{\operatorname{pd}}
\newcommand{\Gpd}{\operatorname{Gpd}}

\newcommand{\Gdim}{\operatorname{G--dim}}

\newcommand{\Ext}{\operatorname{Ext}}
\newcommand{\Supp}{\operatorname{Supp}}

\newcommand{\Hom}{\operatorname{Hom}}

\newcommand{\Ann}{\operatorname{Ann}}

\newcommand{\depth}{\operatorname{depth}}
\newcommand{\width}{\operatorname{width}}

\newcommand{\lo}{\longrightarrow}
\newcommand{\fm}{\frak{m}}
\newcommand{\fp}{\frak{p}}
\newcommand{\fq}{\frak{q}}
\newcommand{\fa}{\frak{a}}

\newenvironment{prf}[1][Proof]{\begin{proof}[\bf #1]}{\end{proof}}

\begin{document}

\author[K. Divaani-Aazar, F. Mohammadi Aghjeh Mashhad and M. Tousi]
{Kamran Divaani-Aazar, Fatemeh Mohammadi Aghjeh Mashhad and Massoud Tousi}

\title[On the existence of certain modules of finite Gorenstein homological dimensions]
{On the existence of certain modules of finite Gorenstein homological dimensions}

\address{K. Divaani-Aazar, Department of Mathematics, Az-Zahra University, Vanak, Post Code
19834, Tehran, Iran-and-School of Mathematics, Institute for Research in Fundamental Sciences
(IPM), P.O. Box 19395-5746, Tehran, Iran.}
\email{kdivaani@ipm.ir}

\address{F. Mohammadi Aghjeh Mashhad,  Parand Branch, Islamic Azad University, Tehran, Iran.}
\email{mohammadi\_fh@yahoo.com}

\address{M. Tousi, Department of Mathematics, Shahid Beheshti University, Tehran, Iran-and-School
of Mathematics, Institute for Research in Fundamental Sciences (IPM), P.O. Box 19395-5746,
Tehran, Iran.} \email{mtousi@ipm.ir}


\keywords {Auslander categories; Cohen-Macaulay complexes; G-dimension; Gorenstein injective
dimension; Gorenstein injective modules; normalized dualizing complexes; totally reflexive modules.\\
The third author was supported by a grant from IPM (No. 90130211).}

\begin{abstract} Let $(R,\fm)$ be a commutative Noetherian local ring. It is known that $R$ is
Cohen-Macaulay if there exists either a nonzero finitely generated $R$-module of finite injective
dimension or a nonzero Cohen-Macaulay $R$-module of finite projective dimension. In this paper, we
investigate the Gorenstein analogues of these facts.
\end{abstract}

\maketitle

\section{Introduction}

Let $(R,\fm)$ be a commutative Noetherian local ring. It is known that the existence of
certain nonzero $R$-modules with finite homological dimensions forces $R$ to be Cohen-Macaulay or
Gorenstein. More precisely, one has
\begin{enumerate}
\item[i)] if there exists a nonzero $R$-module of finite injective dimension, then $R$ is
Cohen-Macaulay,
\item[ii)] if there exists a nonzero Cohen-Macaulay $R$-module of finite projective dimension,
then $R$ is Cohen-Macaulay,
\item[iii)] if there exists a nonzero cyclic $R$-module of finite injective dimension, then $R$
is Gorenstein.
\end{enumerate}
The first assertion is known as Bass's Theorem. The analogue of iii) for Gorenstein homological
dimensions is shown to be true; see \cite[Theorem 4.5]{FF}. The analogues of i) and ii) for Gorenstein
homological dimensions are open questions; see \cite[Pages 40 and 147]{C}, \cite{T2} and
\cite[Questions 1.31 and 3.26]{CFoH}.

\begin{question} If there exists a nonzero finitely generated $R$-module with finite Gorenstein
injective dimension, is then $R$ Cohen-Macaulay?
\end{question}

\begin{question} If there exists a nonzero Cohen-Macaulay $R$-module with finite G-dimension, is
then $R$ Cohen-Macaulay?
\end{question}

Takahashi has given some partial answers to the above questions in \cite{T1} and \cite{T2}. Also, he
showed that if the answer to Question 1.1 is affirmative, then so is the answer to Question 1.2. This
paper has two goals: firstly, it is to show that the assumptions of the above questions are not completely
independent of each other and secondly, it is to give some partial answers to the above questions.

In Section 3, we are concerned with our first goal. When $R$ possesses a normalized dualizing complex,
we show that there exists a nonzero finitely generated $R$-module of finite Gorenstein injective
dimension if and only if there exists a Cohen-Macaulay complex of finite G-dimension; see Theorem
3.7 below. Also, without assuming the existence of a normalized dualizing complex, we prove that
there exists a nonzero Cohen-Macaulay $R$-module of finite Gorenstein injective dimension if and
only if there exists a nonzero Cohen-Macaulay $R$-module of finite G-dimension, see Theorem 3.8 below.
Among other things, this indicates that an affirmative answer to Question 1.1 yields an affirmative
answer to Question 1.2. This was previously shown in \cite{T2}.

In Section 4, we deal with our second goal. If $L$ is a nonzero finitely generated $R$-module
such that either its Gorenstein injective dimension is finite or $L$ is Cohen-Macaulay and its
G-dimension is finite, then we show that $\dim_RL=\depth R-\grade_RL$. Also, if $L$
is a nonzero Cohen-Macaulay $R$-module such that either its Gorenstein injective dimension or
its G-dimension is finite, then we prove that $\dim R/\fp+\depth R_{\fp}=\depth R$ for all
$\fp\in \Ass_RL$. These two facts immediately yield the following partial answers to Questions
1.1 and 1.2: \\Assume that $R$ possesses
a nonzero finitely generated $R$-module $L$ such that either:\\
$L$ has finite Gorenstein injective dimension and $\dim_RL=\dim R-\grade_RL$,\\
$L$ is Cohen-Macaulay with finite G-dimension and $\dim_RL=\dim R-\grade_RL$,\\
$L$ is Cohen-Macaulay with finite Gorenstein injective dimension and
$\dim R/\fp+\depth R_{\fp}=\dim R$ for some $\fp\in \Ass_RL$; or\\
$L$ is Cohen-Macaulay with finite G-dimension and $\dim R/\fp+\depth R_{\fp}=\dim R$ for some
$\fp\in \Ass_RL$.\\
Then $R$ is Cohen-Macaulay.

\section{Prerequisites}

Throughout this paper, $(R,\fm,k)$ is a commutative Noetherian local ring with
nonzero identity. The $\fm$-adic completion of $R$ will be denoted by $\widehat{R}$.

\vspace{0.3cm}
{\bf A. Hyperhomology}
\vspace{0.3cm}

As we will use technical side of hyperhomology and derived category of $R$-modules,
$\mathcal{D}(R)$, we recall some necessary information and notations which are needed
in this paper. The objects in  $\mathcal{D}(R)$ are complexes of $R$-modules and
symbol $\simeq$ denotes isomorphisms in this category. For a complex
$$\textbf{X}=\cdots \lo \textbf{X}_{n+1} \overset{{\partial}_{n+1}^\textbf{X}}\lo
\textbf{X}_n \overset{{\partial}_n^\textbf{X}}\lo \textbf{X}_{n-1}\lo \cdots$$ in
$\mathcal{D}(R)$, its \emph{supremum} and  \emph{infimum} are defined, respectively,
by $\sup \textbf{X}:=\sup \{i\in \mathbb{Z}|H_i(\textbf{X})\neq 0\}$ and $\inf
\textbf{X}:=\inf \{i\in \mathbb{Z}|H_i(\textbf{X})\neq 0\}$, with the usual convention
that $\sup \emptyset=-\infty$ and $\inf \emptyset=\infty$. For an integer $\ell$,
$\Sigma^{\ell}\textbf{X}$ is the complex $\textbf{X}$ shifted $\ell$ degrees to the
left. A complex $\textbf{X}$ is said to be \emph{non-exact} if it has some nonzero
homology modules. The full subcategory of complexes homologically bounded to the right
(resp. left) is denoted by $\mathcal{D}_{\sqsupset}(R)$ (resp. $\mathcal{D}_{\sqsubset}(R)$).
Also, the full subcategories of homologically bounded complexes and of complexes with
finitely generated homology modules will be denoted by $\mathcal{D}_{\Box}(R)$ and
$\mathcal{D}^f(R)$, respectively. We may and do identify the category of $R$-modules and
$R$-homomorphisms with, $\mathcal{D}_0(R)$, the subcategory of complexes homologically
concentrated in degree zero. Finally, we set $\mathcal{D}_{\sqsupset}^f(R):=
\mathcal{D}_{\sqsupset}(R)\cap \mathcal{D}^f(R)$ and  $\mathcal{D}_{\Box}^f(R):=
\mathcal{D}_{\Box}(R)\cap \mathcal{D}^f(R)$. Recall that for any complex $\textbf{X}$,
$\depth_R\textbf{X}:=-\sup{\bf R}\Hom_R(k,\textbf{X})$, $\width_R\textbf{X}:=
\inf(k\otimes^{\bf L}_R\textbf{X})$, $\Supp_R\textbf{X}:=\underset{i\in \mathbb{Z}}\bigcup
\Supp_R(H_i(\textbf{X}))$ and $\dim_R \textbf{X}:=\sup\{\dim R/\fp-\inf
\textbf{X}_{\fp}|\fp\in \Supp_R\textbf{X}\}$.

Let $\textbf{X}\in \mathcal{D}_{\sqsupset}(R)$ and/or $\textbf{Y}\in \mathcal{D}_{\sqsupset}(R)$.
The left-derived tensor product complex of $\textbf{X}$ and $\textbf{Y}$ in $\mathcal{D}(R)$ is
denoted by $\textbf{X}\otimes_R^{{\bf L}}\textbf{Y}$ and is defined by
$$\textbf{X}\otimes_R^{{\bf L}}\textbf{Y}\simeq \textbf{F}\otimes_R\textbf{Y}\simeq \textbf{X}
\otimes_R\textbf{F}^{'}\simeq \textbf{F}\otimes_R\textbf{F},^{'}$$ where $\textbf{F}$ and
$\textbf{F}^{'}$ are flat resolutions of $\textbf{X}$ and $\textbf{Y}$, respectively.
Also, let $\textbf{X}\in \mathcal{D}_{\sqsupset}(R)$ and/or $\textbf{Y}\in \mathcal{D}_
{\sqsubset}(R)$. The right derived homomorphism complex of $\textbf{X}$ and $\textbf{Y}$ in
$\mathcal{D}(R)$ is denoted by ${\bf R}\Hom_R(\textbf{X},\textbf{Y})$ and is defined by $${\bf R}\Hom_R(\textbf{X},\textbf{Y})\simeq \Hom_R(\textbf{P},\textbf{Y})\simeq
\Hom_R(\textbf{X}, \textbf{I})\simeq \Hom_R(\textbf{P},\textbf{I}),$$ where $\textbf{P}$ and
$\textbf{I}$ are projective resolution of $\textbf{X}$ and injective resolution of $\textbf{Y}$,
respectively.

Next, we recall another right derived functor. Let $\fa$ be an ideal of $R$. The right
derived functor of $\fa$-section functor $\Gamma_{\fa}(-)={\varinjlim}_n\Hom_R(R/\fa^n,-)$
is denoted by ${\bf R}\Gamma_{\fa}(-)$. For any complex $\textbf{X}\in
\mathcal{D}_{\sqsubset}(R)$, the complex ${\bf R}\Gamma_{\fa}(\textbf{X})\in
\mathcal{D}_{\sqsubset}(R)$ is defined by ${\bf R}\Gamma_{\fa}(\textbf{X}):=
\Gamma_{\fa}(\textbf{I})$, where $\textbf{I}$ is an (every) injective resolution of
$\textbf{X}$. Also, for any two complexes $\textbf{X}\in \mathcal{D}_{\sqsupset}(R)$ and
$\textbf{Y}\in\mathcal{D}_{\sqsubset}(R)$, ${\bf R}\Gamma_{\fa}(\textbf{X},\textbf{Y})$ is
defined by ${\bf R}\Gamma_{\fa}(\textbf{X},\textbf{Y}):={\bf R}\Gamma_{\fa}({\bf R}
\Hom_R(\textbf{X},\textbf{Y}))$.

\vspace{0.3cm}
{\bf B. Gorenstein Homological Dimensions}
\vspace{0.3cm}

Here, we recall some definitions of Gorenstein homological algebra. A finitely generated
$R$-module $M$ is said to be \emph{totally reflexive} if there exists an exact complex
$\textbf{F}$ of finitely generated free $R$-modules such that $M\cong \im(\textbf{F}_0
\lo \textbf{F}_{-1})$ and $\Hom_R(\textbf{F},R)$ is exact. An $R$-module $N$ is said to
be \emph{Gorenstein injective} if there exists an exact complex $\textbf{I}$ of injective
$R$-modules such that $N\cong \im(\textbf{I}_1\lo \textbf{I}_0)$ and $\Hom_R(E,\textbf{I})$
is exact for all injective $R$-modules $E$. For a complex $\textbf{X}\in \mathcal{D}^{f}_
{\sqsupset}(R)$ (resp. $\textbf{X}\in \mathcal{D}_{\sqsubset}(R)$ ), its \emph{G-dimension}
(resp. \emph{Gorenstein injective dimension}) is defined by $$\begin{array}{llll}\Gdim_R
\textbf{X}&:=\inf\{\sup\{l\in \mathbb{Z}|\textbf{Q}_l\neq 0\}|\textbf{Q}
\  \ \text{is a bounded to the right complex of} \\
&\text{totally reflexive R-modules such that} \   \textbf{Q}\simeq \textbf{X} \},
\end{array}
$$ respectively,
$$\begin{array}{llll}\Gid_R\textbf{X}&:=\inf\{\sup\{l\in \mathbb{Z}|
\textbf{E}_{-l}\neq 0\}|\textbf{E}\  \ \text{is a bounded to the left complex of} \\
&\text{Gorenstein injective R-modules such that} \   \textbf{X}\simeq \textbf{E}\}.
\end{array}
$$

\vspace{0.3cm}
{\bf C. Auslander Categories}
\vspace{0.3cm}

A \emph{normalized dualizing complex} of $R$ is a complex $\textbf{D}\in \mathcal{D}_{\Box}^f(R)$
such that the homothety morphism $R\longrightarrow {\bf R}\Hom_R(\textbf{D},\textbf{D})$
is an isomorphism in $\mathcal{D}(R)$, $\textbf{D}$ has finite injective dimension and
$\sup \textbf{D}=\dim R$. From the definition, it is obvious that
for any normalized dualizing complex $\textbf{D}$, one has $\Supp_R\textbf{D}=\Spec R$.

Assume that $R$ possesses a normalized dualizing complex $\textbf{D}$. The  \emph{Auslander categories}
with respect to $\textbf{D}$, denoted by $\mathcal{A}(R)$ and $\mathcal{B}(R)$, are the full
subcategories of $\mathcal{D}_{\Box}(R)$ whose objects are specified as follows:
\begin{align*}
\mathcal{A}(R)&\,=\, \left\{ \textbf{X} \in \mathcal{D}_{\Box}(R) \: \left|
\begin{array}{l}
\mbox{ $\eta_\textbf{X}:\textbf{X}\longrightarrow {\bf R}\Hom_R(\textbf{D},\textbf{D}\otimes_R^{\bf L}
\textbf{X})$ is an iso-} \\
\mbox{ morphism in $\mathcal{D}(R)$ and $\textbf{D}\otimes_R^{\bf L}\textbf{X }\in
\mathcal{D}_{\Box}(R)$}
\end{array}
\right.
\right\},
\end{align*}
and
\begin{align*}
\mathcal{B}(R) & \,=\, \left\{ \textbf{Y} \in \mathcal{D}_{\Box}(R) \: \left|
\begin{array}{l}
\mbox{ $\varepsilon_\textbf{Y}:\textbf{D}\otimes_R^{\bf L}{{\bf R}\Hom_R(\textbf{D},
\textbf{Y})}\longrightarrow \textbf{Y}$ is an
isomor-} \\
\mbox{ phism in $\mathcal{D}(R)$ and ${\bf R}\Hom_R(\textbf{D},\textbf{Y})\in
\mathcal{D}_{\Box}(R)$}
\end{array}
\right.
\right\}.
\end{align*}

\section{Connections between Questions 1.1 and 1.2}

We begin this section with the following result which relates the finiteness of Gorenstein
homological dimensions to Auslander categories.

\begin{lemma} Let $(R,\fm)$ be a local ring possessing a normalized dualizing complex, $\textbf{X}\in
\mathcal{D}^f_{\Box}(R)$ and $\textbf{Y}\in \mathcal{D}_{\Box}(R)$. The following hold.
\begin{enumerate}
\item[i)]  $\textbf{X}\in \mathcal{A}(R)$ if and only if
$\Gdim_R\textbf{X}<\infty$; and
\item[ii)]  $\textbf{Y}\in \mathcal{B}(R)$ if and only if $\Gid_R\textbf{Y}
<\infty$.
\end{enumerate}
\end{lemma}

\begin{prf} i) See \cite[Proposition 3.8 b) and Theorem 4.1]{CFrH}.

ii) See \cite[Theorem 4.4]{CFrH}.
\end{prf}

The next result computes $-\inf {\bf R}\Gamma_{\fm}(\textbf{X},\textbf{Y})$ and $-\sup {\bf R}\Gamma_{\fm}(\textbf{X},\textbf{Y})$ in some special cases.

\begin{lemma}  Let $(R,\fm)$ be a local ring possessing a normalized dualizing complex $\textbf{D}$
and $\textbf{X},\textbf{Y}\in \mathcal{D}_{\Box}^f(R)$.
\begin{enumerate}
\item[i)] If either projective dimension of $\textbf{X}$
or $\textbf{Y}$ is finite, then $$-\inf {\bf R}\Gamma_{\fm}(\textbf{X},\textbf{Y})=\sup {\bf R}\Hom_R(\textbf{Y},\textbf{D}\otimes_R^{{\bf L}}\textbf{X})$$ and $$-\sup {\bf R}\Gamma_{\fm}(\textbf{X},\textbf{Y})=\inf {\bf R}\Hom_R(\textbf{Y},\textbf{D}
\otimes_R^{{\bf L}}\textbf{X}).$$
\item[ii)] If projective dimension of $\textbf{X}$ and Gorenstein injective dimension of $\textbf{Y}$
are finite, then $$-\inf {\bf R}\Gamma_{\fm}(\textbf{X},\textbf{Y})=\sup {\bf R}\Hom_R({\bf R}\Hom_R(\textbf{D},\textbf{Y}),\textbf{X})$$ and $$-\sup {\bf R}\Gamma_{\fm}(\textbf{X},\textbf{Y})=\inf {\bf R}\Hom_R({\bf R}\Hom_R(\textbf{D},\textbf{Y}),\textbf{X}).$$
\end{enumerate}
\end{lemma}

\begin{prf} See \cite[Corollary 3.4]{MD}.
\end{prf}

\begin{lemma} Let $(R,\fm)$ be a local ring possessing a normalized dualizing complex
$\textbf{D}$ and $\textbf{X}\in \mathcal{D}_{\Box}^f(R)$ a non-exact complex. Then $\dim_RX=\sup {\bf R}\Hom_R(\textbf{X},\textbf{D})$ and $\depth_RX=\inf {\bf R}\Hom_R(\textbf{X},\textbf{D})$. In particular,
$\depth R=\inf \textbf{D}$.
\end{lemma}

\begin{prf} See \cite[Theorem 15.10 and 16.20 a) and b)]{Fo1}.
\end{prf}

The following lemma states some properties of modules with finite G-dimension over a local ring
with a normalized dualizing complex.

\begin{lemma} Let $(R,\fm)$ be a local ring possessing a normalized dualizing complex $\textbf{D}$
and $M$ a nonzero finitely generated $R$-module. The following hold.
\begin{enumerate}
\item[i)] $\inf (\textbf{D}\otimes_R^{\bf L}M)=\depth_R {\bf R}\Hom_R(M,R)=\depth R$.\\
Assume that $M$ has finite G-dimension. Then
\item[ii)]  $\depth R\leq \sup (\textbf{D}\otimes_R^{\bf L}M)=\dim_R {\bf R}\Hom_R(M,R)\leq \dim R$.
\item[iii)] $\sup (\textbf{D}\otimes_R^{{\bf L}}M)= \sup \{\dim R/\fp+\Gdim_{R_{\fp}}M_{\fp}|\fp\in
\Supp_RM \}$.
\end{enumerate}
\end{lemma}

\begin{prf} i) \cite[Corollary A.4.16]{C} and Lemma 3.3 yield that $$\inf (\textbf{D}
\otimes_R^{\bf L}M)=\inf \textbf{D}+\inf M=\depth R.$$ By \cite[Proposition 3.8]{Fo3} for
any complex $\textbf{X}\in \mathcal{D}_{\sqsubset}(R)$, one has $-\sup {\bf R}
\Gamma_{\fm}(\textbf{X})=\depth_R\textbf{X}$. Hence,  Lemma 3.2 i) implies that
$$\depth_R{\bf R}\Hom_R(M,R)=-\sup {\bf R}\Gamma_{\fm}(M,R)=
\inf {\bf R}\Hom_R(R,\textbf{D}\otimes_R^{\bf L}M)=\inf (\textbf{D}\otimes_R^{\bf L}M).$$

Now, assume that $M$ has finite G-dimension.

ii) Lemma 3.1 i) yields that $M\in \mathcal{A}(R)$ and then by \cite[Proposition 3.3.7 a)]{C}
and Lemma 3.3, one has $$\depth R=\inf \textbf{D}\leq  \sup (\textbf{D}\otimes_R^{\bf L}M)\leq \sup
\textbf{D}=\dim R.$$ Also, \cite[Theorem 2.2.3 and Definition 2.1.6]{C} yield that ${\bf R}
\Hom_R(M,R)\in \mathcal{D}_{\Box}^f(R)$. On the other hand, \cite[Proposition 3.14 d)]{Fo3}
implies that for any complex $\textbf{X}\in \mathcal{D}^f_{\Box}(R)$, one has
$-\inf {\bf R}\Gamma_{\fm}(\textbf{X})=\dim_R\textbf{X}$. Hence Lemma 3.2 i) asserts that
$$\dim_R{\bf R}\Hom_R(M,R)=-\inf {\bf R}\Gamma_{\fm}(M,R)=\sup {\bf R}\Hom_R(R,\textbf{D}
\otimes_R^{\bf L}M)=\sup (\textbf{D}\otimes_R^{\bf L}M).$$

iii) Set $N:=R$ in the proof of \cite[Lemma 4.2 i)]{MD}.
\end{prf}

Assume a finitely generated $R$-module $N$ with finite injective dimension is given.
The problem of finding a finitely generated $R$-module with finite projective
dimension and the same support as $N$ has been investigated in \cite[Theorem 5.7]{PS}
and \cite{Sh}. The following is Gorenstein analogue of this problem. It will be needed
in the proof of our main results.

\begin{lemma} Let $(R,\fm)$ be a local ring with a normalized dualizing complex $\textbf{D}$.
Assume that $N$ is a nonzero finitely generated $R$-module with finite Gorenstein injective
dimension. Set $H:=H_{-\depth R}({\bf R}\Hom_R(\textbf{D},N))$. Then
\begin{enumerate}
\item[i)] $\inf {\bf R}\Hom_R(\textbf{D},N)=-\depth R=\sup {\bf R}\Hom_R(\textbf{D},N)$,
\item[ii)] $H$ is a nonzero finitely generated $R$-module and $H\simeq
\Sigma^{\depth R}{\bf R}\Hom_R(\textbf{D},N)$,
\item[iii)] $\Gdim_RH<\infty$; and
\item[iv)] $\depth_RH=\depth_RN$ and $\Supp_RH=\Supp_RN$.
\end{enumerate}
\end{lemma}

\begin{prf} i) Lemma 3.1 ii) yields that $N\in \mathcal{B}(R)$, and so
$$N\simeq D\otimes_R^{{\bf L}}{\bf R}\Hom_R(\textbf{D},N) \   \  (*),$$ and
${\bf R}\Hom_R(\textbf{D},N)\in \mathcal{D}_{\Box}(R)$.  Also, by \cite[Lemma A.4.4]{C}
all homology modules of ${\bf R}\Hom_R(\textbf{D},N)$ are finitely
generated. Applying \cite[Corollary A.4.16]{C} and Lemma 3.3 to $(*)$, implies that
$$\inf {\bf R}\Hom_R(\textbf{D},N)=-\inf \textbf{D}=-\depth R.$$  Also, \cite[ A.4.6.1]{C}
implies that $$\sup {\bf R}\Hom_R(\textbf{D},N)\leq \sup N-\inf \textbf{D}=-\depth R.$$ Hence,
$$\inf {\bf R}\Hom_R(\textbf{D},N)=-\depth R=\sup {\bf R}\Hom_R(\textbf{D},N).$$

ii) is clear by i).

iii)  Since $N\in \mathcal{B}(R)$, \cite[Theorem 3.3.2 b)]{C} yields that ${\bf R}
\Hom_R(\textbf{D},N)\in \mathcal{A}(R)$, and so $H\in \mathcal{A}(R)$. Thus by Lemma 3.1 i),
it turns out that $\Gdim_RH<\infty$.

iv) For any complex $\textbf{X}$ and any integer $s$, it is easy to see that
$\depth_R\Sigma^{s}\textbf{X}=\depth_R\textbf{X}-s$. So, by \cite[Lemma A.6.4 and A.6.3.2]{C}
and Lemma 3.3, we have
\begin{align*}
\depth_RH&=\depth_R({\bf R}\Hom_R(\textbf{D},N))-\depth R \\
&=\width_R\textbf{D}+\depth_RN-\depth R\\
&=\inf \textbf{D}+\depth_RN-\depth R\\
&=\depth_RN.
\end{align*}
By applying \cite[Page 36]{Fo2} to $(*)$ and using ii), one has
$$\begin{array}{ll}
\Supp_RN&=\Supp_R\textbf{D}\cap \Supp_R({\bf R}\Hom_R(\textbf{D},N))\\
&=\Spec R\cap \Supp_R({\bf R}\Hom_R(\textbf{D},N))\\
&=\Supp_RH.
\end{array}$$
Recall that $\Supp_R\textbf{D}=\Spec R$.
\end{prf}

\begin{definition} Let $(R,\fm)$ be a local ring. A non-exact complex $\textbf{X}\in
\mathcal{D}_{\Box}^f(R)$ is said to be Cohen-Macaulay if $\depth_R\textbf{X}=\dim_R
\textbf{X}$.
\end{definition}

Now, we are ready to present the main results of this section, which explain the relationship between
the assumptions of Questions $1.1$ and $1.2$.

\begin{theorem} Let $(R,\fm)$ be a local ring possessing a normalized dualizing complex
$\textbf{D}$. The following are equivalent:
\begin{enumerate}
\item[i)] There exists a nonzero finitely generated $R$-module $N$ such that $\Gid_RN<\infty$.
\item[ii)] There exists a nonzero finitely generated $R$-module $M$ such that $\Gdim_RM< \infty$
and ${\bf R}\Hom_R(M,R)$ is a Cohen-Macaulay complex.
\item[iii)] There exists a Cohen-Macaulay complex $\textbf{X}$ such that $\Gdim_R\textbf{X}< \infty$.
\end{enumerate}
\end{theorem}

\begin{prf} $i\Rightarrow ii)$ Set $s:=\depth R$ and $M:=H_{-s}({\bf R}\Hom_R(\textbf{D},N))$.
By Lemma 3.5, $M$ is a nonzero finitely generated $R$-module, $M\simeq \Sigma^s{\bf R}
\Hom_R(\textbf{D},N)$ and $\Gdim_RM<\infty$. Lemma 3.1 ii) implies that $N\in \mathcal{B}(R)$,
and so $N\simeq \textbf{D}\otimes_R^{{\bf L}}{\bf R}\Hom_R(\textbf{D},N)$. Hence,
\begin{align*}
\sup (\textbf{D}\otimes_R^{\bf L}M) &= \sup (\textbf{D}\otimes_R^{\bf L}\Sigma^s{\bf R}
\Hom_R(\textbf{D},N))\\
&\overset{(a)}=\sup \Sigma^s(\textbf{D}\otimes_R^{\bf L}{\bf R}\Hom_R(\textbf{D},N)) \\
&=\sup \Sigma^sN\\
&=s,
\end{align*}
where $(a)$ is due to \cite[A.2.4.3]{C}. By \cite[Theorem 2.2.3 and Definition 2.1.6]{C}, we have
${\bf R}\Hom_R(M,R)\in \mathcal{D}_{\Box}^f(R)$, and so by Lemma 3.4 i) and ii) the complex
${\bf R}\Hom_R(M,R)$ is Cohen-Macaulay.

$ii\Rightarrow iii)$ Set $\textbf{X}:={\bf R}\Hom_R(M,R)$. \cite[Theorem 2.1.10 a) and Corollary 2.3.8]{C}
yields that $\Gdim_R\textbf{X}<\infty$. Since by the assumption $\textbf{X}$ is Cohen-Macaulay,
iii) follows.

$iii\Rightarrow i)$ Since $\Gdim_R\textbf{X}<\infty$, Lemma 3.1 i) implies that
$\textbf{X}\in \mathcal{A}(R)$, and so ${\bf R}\Hom_R(\textbf{X},\textbf{D})\in \mathcal{B}(R)$ by
\cite[Proposition 4.1 v)]{CH}. Hence, $\Gid_R({{\bf R}\Hom_R(\textbf{X},\textbf{D})})<\infty$ by
Lemma 3.1 ii).

As $\textbf{X}$ is a Cohen-Macaulay complex, Lemma 3.3 yields that $\sup{\bf R}\Hom_R(\textbf{X},
\textbf{D})=\inf{\bf R}\Hom_R(\textbf{X},\textbf{D})$. Set $d:=\sup{\bf R}\Hom_R(\textbf{X},
\textbf{D})$ and $N:=H_{d}({\bf R}\Hom_R(\textbf{X},\textbf{D}))$. Then $N\simeq \Sigma^{-d}{\bf R}\Hom_R(\textbf{X},\textbf{D})$.  Therefore, $N$ is a nonzero finitely
generated $R$-module of finite Gorenstein injective dimension.
\end{prf}

\begin{theorem} Let $(R,\fm)$ be a local ring. The following are equivalent:
\begin{enumerate}
\item[i)] There exists a nonzero Cohen-Macaulay $R$-module with finite G-dimension.
\item[ii)] There exists a nonzero finite length $R$-module with finite G-dimension.
\item[iii)] There exists a nonzero finite length $R$-module with finite Gorenstein injective dimension.
\item[iv)] There exists a nonzero Cohen-Macaulay $R$-module with finite Gorenstein injective dimension.
\end{enumerate}
\end{theorem}

\begin{prf} $ii\Rightarrow i)$ and $iii\Rightarrow iv)$ are trivial.

$i\Rightarrow ii)$ Let $M$ be a nonzero Cohen-Macaulay $R$-module with finite G-dimension. Set
$t:=\depth_RM$ and assume that $x_1,\ldots, x_t$ is an $M$-regular sequence. Then $M/{(x_1,\ldots x_t)}M$
is a finite length $R$-module, and also by \cite[Theorem 8.7 7)]{A},  $\Gdim_R{M/{(x_1,\ldots x_t)}M}
=\Gdim_RM+t$. This completes the proof.

$ii\Longleftrightarrow iii)$ Let $A$ be an Artinian $R$-module. Then $A$ possesses the structure of
an $\widehat{R}$-module in a natural way such that $\ell_R(A)=\ell_{\widehat{R}}(A)$. On the other
hand, any Artinian $\widehat{R}$-module $B$ is also Artinian when it is considered as an $R$-module
via the natural ring homomorphism $R\lo \widehat{R}$ and $\ell_{\widehat{R}}(B)=\ell_R(B)$. For any
finite length $R$-module $C$, one has $C\cong C\otimes_R\widehat{R}$, and so \cite[Theorem 8.7 5)]{A}
implies that $\Gdim_RC=\Gdim_{\widehat{R}}C$. Also, for any Artinian $R$-module $C$, \cite[Lemma 3.6]{Sa}
yields that $\Gid_RC=\Gid_{\widehat{R}}C$. Hence, in view of \cite[Proposition 3.4.9]{St}, we may and
do assume that $R$ is complete. Let $(-)^{\vee}:=\Hom_R(-,E_R(R/\fm))$ and $C$ be a nonzero finite
length $R$-module. Then $C\cong (C^{\vee})^{\vee}$ and $\ell_R(C^{\vee})<\infty$. Thus by
\cite[Proposition 4.1 v) and vi)]{CH}, it turns out that $C\in \mathcal{A}(R)$ (resp. $C\in \mathcal{B}(R)$)
if and only if $C^{\vee}\in \mathcal{B}(R)$ (resp. $C^{\vee}\in \mathcal{A}(R)$). Now, Lemma 3.1
completes the argument of this part.

$iv\Rightarrow iii)$ Let $M$ be a nonzero Cohen-Macaulay $R$-module with finite Gorenstein injective
dimension. Let $n:=\dim_RM$ and $\underline{x}:=x_1, x_2,\ldots, x_n\in \fm$ be a maximal $M$-regular
sequence. As $\ell_R(M/(\underline{x})M)<\infty$, it is enough to show that $\Gid_R(M/(\underline{x})M)<\infty$.
We do this by using induction on $n$. For $n=0$, there is nothing to prove. Now, assume that $n>0$ and
the assertion holds for $n-1$. Set $N:=M/x_1M$. Then $N$ is a nonzero Cohen-Macaulay $R$-module with
$\dim_RN=n-1$. By applying \cite[Theorem 2.25]{H} on the short exact sequence $$0\lo M\overset{x_1}\lo M\lo
N\lo 0,$$ we deduce that $\Gid_RN<\infty$. Since $x_2, x_3, \ldots, x_n$ is a maximal $N$-regular sequence
and $$N/(x_2, x_3, \ldots, x_n)N\cong M/(\underline{x})M,$$ from the induction hypothesis, we conclude that $\Gid_R(M/(\underline{x})M)<\infty$, as desired.
\end{prf}

\section{Some partial answers to Questions 1.1 and 1.2}

In this section, we give some partial answers to Questions 1.1 and 1.2. We begin this section
with the following lemma which is the dual of Lemma 3.5.

\begin{lemma} Let $(R,\fm)$ be a local ring with a normalized dualizing complex $\textbf{D}$.
Assume that $M$ is a nonzero Cohen-Macaulay $R$-module with finite G-dimension. Set $H:=H_{\depth R}(\textbf{D}\otimes_R^{{\bf L}}M)$. Then
\begin{enumerate}
\item[i)] $\sup (\textbf{D}\otimes_R^{\bf L}M)=\depth R=\inf (\textbf{D}\otimes_R^{\bf L}M)$,
\item[ii)] $H$ is a nonzero finitely generated $R$-module and $H\simeq \Sigma^{-\depth R}(\textbf{D}
\otimes_R^{{\bf L}}M)$,
\item[iii)] $\Gid_RH<\infty$; and
\item[iv)] $\depth_RH=\depth_RM$ and $\Supp_RH=\Supp_RM$, and so $H$ is Cohen-Macaulay.
\end{enumerate}
\end{lemma}

\begin{prf} i) Lemma 3.4 i) yields that $\inf (\textbf{D}\otimes_R^{\bf L}M)=\depth R$. On the other hand,
we have
\begin{align*}
\depth R&\overset{(a)}\leq  \sup (\textbf{D}\otimes_R^{\bf L}M) \\
&=\sup {\bf R}\Hom_R(R,\textbf{D}\otimes_R^{\bf L}M) \\
&\overset{(b)}=-\inf {\bf R}\Gamma_{\fm}(M,R)\\
&\overset{(c)}\leq -\inf {\bf R}\Hom_R(M,R)+\dim_RM\\
&\overset{(d)}= \Gdim_RM+\dim_RM\\
&=\Gdim_RM+\depth_RM\\
&\overset{(e)}=\depth R,
\end{align*}
where $(a)$ is by Lemma 3.4 ii), $(b)$ is by Lemma 3.2 i) and $(c)$ follows from \cite[Definition 2.1
and Corollary 3.2 ii)]{DH}. $(d)$ holds because by \cite[Theorem 2.2.3]{C}, one has $-\inf {\bf R}\Hom_R(M,R)=\Gdim_RM$, and finally $(e)$ is due to \cite[Theorem 1.4.8]{C}.

ii) Clearly, i) yields that $H$ is nonzero and $H\simeq \Sigma^{-\depth R}(\textbf{D}
\otimes_R^{{\bf L}}M)$. On the other hand,  \cite[Lemma A.4.13]{C} implies that
all homology modules of $\textbf{D}\otimes_R^{{\bf L}}M$ are finitely generated.

iii) Lemma 3.1 i) implies that $M\in \mathcal{A}(R)$. Hence, \cite[Theorem 3.3.2 a)]{C} yields
that $\textbf{D}\otimes_R^{\bf L}M\in \mathcal{B}(R)$, and so $H\in \mathcal{B}(R)$. Now,
Lemma 3.1 ii) implies that $\Gid_RH$ is finite.

iv) Since $M\in \mathcal{A}(R)$, one has $M\simeq {\bf R}\Hom_{R}(\textbf{D},
\textbf{D}\otimes^{\bf L}_{R}M)$, and so
\begin{align*}
\depth_RM&= \depth_R({\bf R}\Hom_R(\textbf{D},\textbf{D}\otimes_R^{{\bf L}}M)) \\
&\overset{(a)}=\width_R\textbf{D}+\depth_R(\textbf{D}\otimes_R^{{\bf L}}M)\\
&\overset{(b)}=\inf \textbf{D}+\depth_R(\textbf{D}\otimes_R^{{\bf L}}M)\\
&\overset{(c)}=\depth R+\depth_R(\textbf{D}\otimes_R^{{\bf L}}M) \\
&\overset{(d)}=\depth_R\Sigma^{-\depth R}(\textbf{D}\otimes_R^{{\bf L}}M)\\
&=\depth_RH,
\end{align*}
where $(a)$ and $(b)$ are due to \cite[Lemma A.6.4 and A.6.3.2]{C} and $(c)$ is by Lemma 3.3.
As for any complex $\textbf{X}$ and any integer $s$, $\depth_R\Sigma^{s}\textbf{X}=
\depth_R\textbf{X}-s$, the equality $(d)$ holds.

Since $\Supp_R\textbf{D}=\Spec R$, by ii) and \cite[Page 36]{Fo2}, one has
$$\begin{array}{ll}
\Supp_RH&=\Supp_R(\textbf{D}\otimes_R^{{\bf L}}M)\\
&=\Supp_R\textbf{D}\cap \Supp_RM\\
&=\Spec R\cap \Supp_RM\\
&=\Supp_RM.
\end{array}$$
\end{prf}

Let $R$ be a local ring and $M$ a nonzero Cohen-Macaulay $R$-module with finite G-dimension. The
first case of the following result can be rephrased as saying that  $M$ is G-perfect in the sense
of Golod. We recall that for any finitely generated $R$-module $L$, $\grade_RL:=\inf\{i\geq 0|
\Ext^{i}_R(L,R)\neq 0\}.$

\begin{theorem} Let $(R,\fm)$ be a local ring and $M$ a nonzero finitely generated $R$-module.
If either $M$ is Cohen-Macaulay with finite G-dimension or $M$ has finite Gorenstein injective
dimension, then $$\dim_RM=\depth R-\grade_RM.$$
\end{theorem}

\begin{prf} Let $T$ be a finitely generated $R$-module. By \cite[Theorem 8.7 5)]{A}, one has $\Gdim_RT=\Gdim_{\widehat{R}}(T\otimes_R{\widehat{R}})$ and if $\Gid_RT<\infty$, then
\cite[Theorem 2.5]{KTY} implies that $\Gid_RT=\Gid_{\widehat{R}}(T\otimes_R\widehat{R})$.
Also, it is easy to see that  $$\grade_RT=\grade_{\widehat{R}}(T\otimes_R{\widehat{R}}).$$ Hence,
without loss of generality, we may and do assume that $R$ is complete, and so $R$ possesses a
normalized dualizing complex $\textbf{D}$. If $M$ is Cohen-Macaulay with finite G-dimension,
then Lemma 4.1 yields that there exists a nonzero finitely generated $R$-module $H$ such that
$\Gid_RH<\infty$ and $\sqrt{\Ann_RH}=\sqrt{\Ann_RM}$. Thus, it suffices to assume that $M$ is
a nonzero finitely generated $R$-module with finite Gorenstein injective dimension.

By Grothendieck's non-vanishing Theorem \cite[Proposition 3.14 d)]{Fo3} and Lemma 3.2 ii),
one has
$$\dim_RM=-\inf {\bf R}\Gamma_{\fm}(M)=-\inf {\bf R}\Gamma_{\fm}(R,M)=\sup {\bf R}\Hom_R({\bf R}\Hom_R(\textbf{D},M),R).$$
But, by \cite[Lemma 4.2 ii)]{MD} we have $$\sup {\bf R}\Hom_R({\bf R}\Hom_R(\textbf{D},M),R)=
\depth R-\grade_RM,$$ which completes the proof.
\end{prf}

\begin{theorem} Let $(R,\fm)$ be a local ring and $M$ a nonzero Cohen-Macaulay $R$-module such that
either its Gorenstein injective dimension or its G-dimension is finite. Then
$$\dim R/\fp+\depth R_{\fp}=\depth R$$ for all $\fp\in \Ass_RM$.
\end{theorem}

\begin{prf}Let $L$ be a nonzero finitely generated $R$-module. Then $L$ is Cohen-Macaulay if and
only if the $\widehat{R}$-module $L\otimes_R\widehat{R}$ is Cohen-Macaulay. \cite[Theorem 8.7 5)]{A}
implies that $\Gdim_RL=\Gdim_{\widehat{R}}(L\otimes_R\widehat{R})$. Also, by \cite[Theorem 2.5]{KTY}
if $\Gid_RL<\infty$, then $\Gid_RL=\Gid_{\widehat{R}}(L\otimes_R\widehat{R})$. Obviously,
$\depth R=\depth \widehat{R}$, and also it is known that
$$\Ass_RL=\{\fq\cap R|\fq\in \Ass_{\widehat{R}}(L\otimes_R\widehat{R})\}.$$ Let $\fp\in \Ass_RL$
and $\fq\in \Ass_{\widehat{R}}(L\otimes_R\widehat{R})$ be such that $\fp=\fq\cap R$. Assume that $L$ is
Cohen-Macaulay. Then, by \cite[Theorem 2.1.2 a)]{BH}, one has $$\dim R/\fp=\dim_RL=\dim_{\widehat{R}}(L\otimes_R{\widehat{R}})=\dim {\widehat{R}}/\fq.$$
Since $\fp{\widehat{R}}\subseteq \fq$ and
$\dim {\widehat{R}}/\fq=\dim {\widehat{R}}/\fp{\widehat{R}}$, one has $\fq\in \Min{(\fp{\widehat{R}})}$.
Then $\depth_{{\widehat{R}}_{\fq}}{\widehat{R}}_{\fq}/\fp{\widehat{R}}_{\fq}=0$, and so by applying
\cite[Proposition 1.2.16 a)]{BH} on the flat local ring extension $R_{\fp}\lo {\widehat{R}}_\fq$,
we have $\depth R_{\fp}=\depth {\widehat{R}}_{\fq}$.
Hence, $$\dim R/\fp+\depth R_{\fp}=\dim {\widehat{R}}/\fq+\depth {\widehat{R}}_{\fq}.$$
So, without loss of generality, we may and do
assume that $R$ is complete. Then $R$ has a normalized dualizing complex $\textbf{D}$.

Suppose that $M$ has finite Gorenstein injective dimension. By Lemma 3.5, there exists a nonzero
Cohen-Macaulay $R$-module $H$ of finite G-dimension such that $\Supp_RH=\Supp_RM$. Then the set
of minimal elements of $\Supp_RH$ and $\Supp_RM$ are equal too. So, $\Ass_RH=\Ass_RM$. Note that
Cohen-Macaulayness of $H$ and $M$ implies that $H$ and $M$ do not have any embedded associated
prime ideals. Hence, we only need to consider the case that $\Gdim_RM$ is finite.
Let $\fp\in \Ass_RM$. Then, by \cite[Proposition 1.3.2 and Theorem 1.4.8]{C}, it turns out that
$$\Gdim_{R_{\fp}}M_{\fp}=\depth R_{\fp}-\depth_{R_{\fp}}M_{\fp}=\depth R_{\fp}.$$ Thus one has:
\begin{align*}
\depth R&\overset{(a)}= \dim_RM + \grade_RM\\
&\leq \dim R/\fp + \grade_{R_{\fp}}M_{\fp}\\
&\leq \dim R/\fp+\depth R_{\fp}\\
&=\dim R/\fp +\Gdim_{R_{\fp}}M_{\fp}\\
&\overset{(b)}\leq \sup (\textbf{D}\otimes_R^{\bf L}M) \\
&\overset{(c)}=\depth R,
\end{align*}
where $(a)$ is by Theorem 4.2, $(b)$ is by Lemma 3.4 iii) and $(c)$ is by Lemma 4.1 i).
\end{prf}

Let $(R,\fm)$ be a local ring and $N$ a nonzero finitely generated $R$-module with finite Gorenstein
injective dimension. It is known that if $\dim N=\dim R$, then $R$ is Cohen-Macaulay; see
\cite[Theorem 1.3]{Y}. The implication $iii)\Rightarrow i)$ in the next result improves this fact.

\begin{corollary} Let $(R,\fm)$ be a local ring. The following are equivalent:
\begin{enumerate}
\item[i)] $R$ is Cohen-Macaulay.
\item[ii)] There exists a nonzero Cohen-Macaulay $R$-module $M$ with finite G-dimension such that
$\dim_RM=\dim R-\grade_RM$.
\item[iii)] There exists a nonzero finitely generated $R$-module $N$ with finite
Gorenstein injective dimension such that  $\dim_RN=\dim R-\grade_RN$.
\item[iv)] There exists a nonzero Cohen-Macaulay $R$-module $M$ of finite G-dimension such that
$\dim R/\fp+\depth R_{\fp}=\dim R$ for some $\fp\in \Ass_RM$.
\item[v)] There exists a nonzero Cohen-Macaulay $R$-module $N$ of finite Gorenstein injective
dimension such that $\dim R/\fp+\depth R_{\fp}=\dim R$ for some $\fp\in \Ass_RN$.
\end{enumerate}
\end{corollary}

\begin{prf} $i\Rightarrow ii)$ and $i\Rightarrow iv)$ are trivial by setting $M:=R$.

$ii\Rightarrow i)$ and $iii\Rightarrow i)$ are clear by Theorem 4.2.

$iv)\Rightarrow i)$ and $v)\Rightarrow i)$ are clear by Theorem 4.3.

$i\Rightarrow iii)$ Set $N:=\Hom_R(R/\textbf{x}R,E(R/\fm))$, where $\textbf{x}=
x_{1},\ldots, x_{n}$ is a system of parameters of $R$. Since $R$ is Cohen-Macaulay,
$\textbf{x}$ is an $R$-regular sequence, and so $\pd_R{R/\textbf{x}R}$ is finite.
Then $\id_RN$ is finite, and so $\Gid_RN=\id_RN$, by \cite[Proposition 3.10]{CFoH}.
Because $N$ has finite length, we have $\dim_RN=0$ and $\grade_RN=\depth R$,
and so the proof is complete.

$i\Rightarrow v)$  Let $N$ be as in the proof $i\Rightarrow iii)$. Then $\Gid_RN<\infty$ and
$N$ has finite length. Thus the assertion follows, by taking $\fp:=\fm$.
\end{prf}

Recall that an ideal of $R$ is said to be {\em complete intersection} if it is generated by
a regular $R$-sequence.

\begin{corollary} Let $(R,\fm)$ be a local ring.
\begin{enumerate}
\item[i)] If there exists a nonzero Cohen-Macaulay $R$-module $M$ with finite G-dimension
such that $\Ann_RM$ is complete intersection, then $R$ is Cohen-Macaulay.
\item[ii)] If there exists a nonzero finitely generated $R$-module $N$ with finite
Gorenstein injective dimension such that $\Ann_RN$ is complete intersection, then $R$ is
Cohen-Macaulay.
\end{enumerate}
\end{corollary}

\begin{prf} i) By the assumption, the ideal $\Ann_RM$ is generated by a regular $R$-sequence
$x_1, x_2,\ldots, x_n$ 'say. Set $T:=R/\Ann_RM$. Then $\depth_TM=\depth_RM$ and $\dim_TM=\dim_RM$,
and so $M$ is also Cohen-Macaulay as a $T$-module. On the other hand, \cite[Theorem 2.2.8]{C}
yields that $$\Gdim_TM=\Gdim_RM-n<\infty.$$ Now, by the implication $ii\Rightarrow i)$ in
Corollary 4.4, it follows that the ring $T$ is Cohen-Macaulay. Hence, $R$ itself is
Cohen-Macaulay.

ii) Let $x_1, x_2,\ldots, x_n$ be a regular $R$-sequence generating the ideal $\Ann_RN$ and set
$T:=R/\Ann_RN$. Then \cite[Theorem 4.2]{BM} implies that $$\Gid_TN=\Gid_RN-n<\infty.$$
Thus, by the implication $iii\Rightarrow i)$ in Corollary 4.4, it turns out that the
ring $T$ is Cohen-Macaulay, and so $R$ is Cohen-Macaulay too.
\end{prf}

\begin{acknowledgement}We would like to thank Henrik Holm and Sean Sather-Wagstaff for
their valuable comments.
\end{acknowledgement}


\end{document}